\documentclass[11pt]{amsart}
\usepackage{dsfont}
\usepackage{times}
\usepackage{enumerate}
\usepackage[euler-digits]{eulervm}
\theoremstyle{plain}
\usepackage{graphicx,color} 
\usepackage{mathtools}

\usepackage{hyperref}
\usepackage{subcaption}
\usepackage{amsmath, amssymb, graphics}
\usepackage{float}

\newtheorem{theorem}{Theorem}[section]
\newtheorem*{theorem*}{Theorem} 

\newtheorem{prop}{Proposition}[section]
\newtheorem{lemma}{Lemma}[section]

\newtheorem{remark}{Remark}[section]

\newtheorem*{conjecture*}{Conjecture}

\begin{document}

\title{Symmetry Breaking Bifurcation of Membranes with Boundary}
\author{Bennett Palmer and \'Alvaro P\'ampano}
\date{\today}
\maketitle

\begin{abstract} We use a bifurcation theory due to Crandall and Rabinowitz to  show the existence of a symmetry breaking bifurcation of a specific one parameter family of axially symmetric disc type solutions of a membrane equation with fixed boundary.  In place  of working directly with the fourth order membrane equation, it is replaced by a second order reduction found in \cite{PP2}.
\\

\noindent K{\tiny EY} W{\tiny ORDS}.\, Bifurcation, Helfrich Energy, Symmetry Breaking.

\noindent MSC C{\tiny LASSIFICATION} (2020).\, 35B32, 49Q10, 53A05.
\end{abstract}

\begin{figure}[H]
\begin{center}
\includegraphics[width=0.75\linewidth]{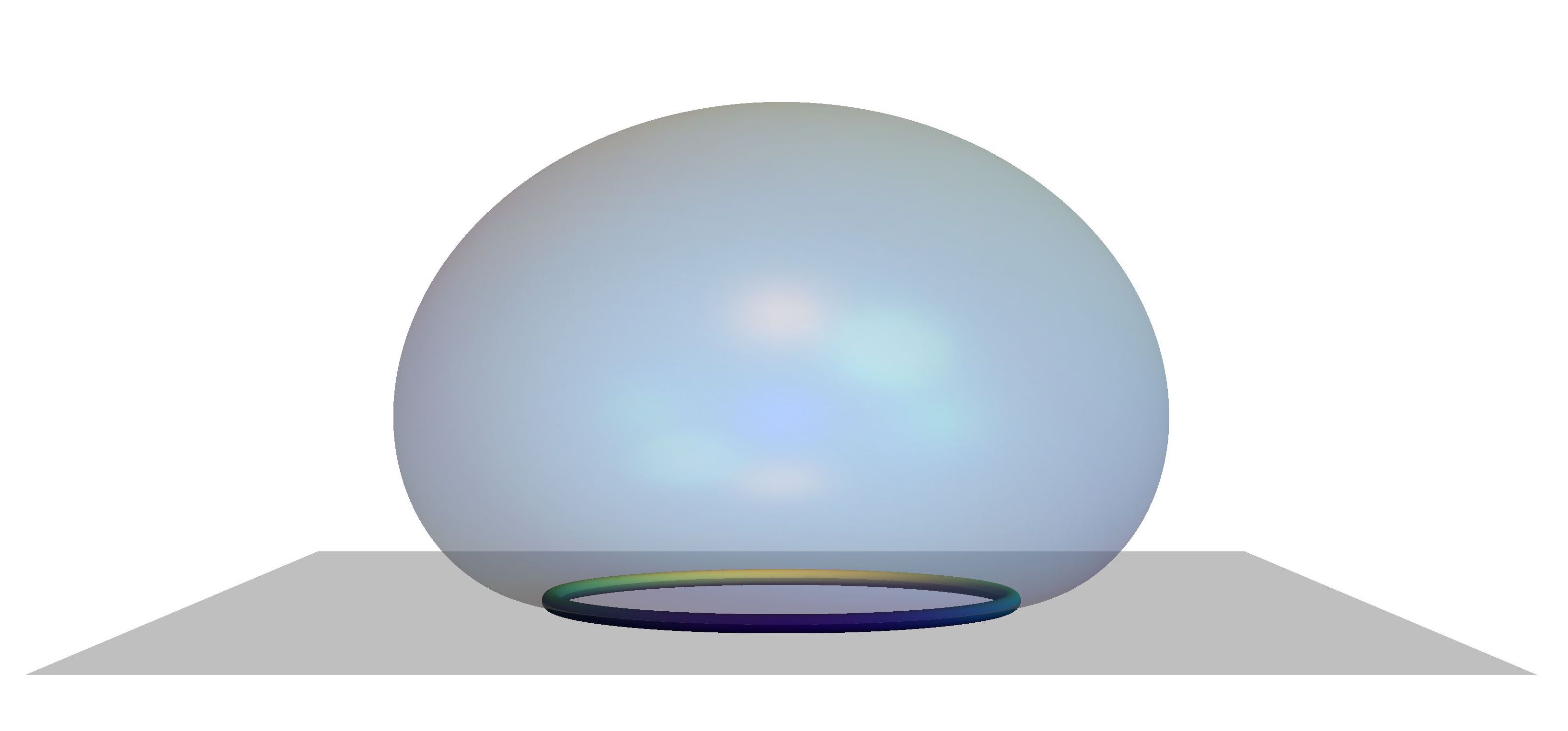}
\end{center} 
\end{figure}

\section{Introduction}

Bilipid membranes are a fundamental structure in cell biology. These membranes are composed of bipolar molecules having a hydrophilic head and a hydrophobic tail. When a sufficiently high density of these molecules is reached in an aqueous solution, the molecules self-assemble so as to hide the tails from the surrounding liquid. This results in the formation of a membrane made of two layers of molecules.

Historically, there has been much interest in explaining the morphology of these membranes. Since closed membranes could form  non spherical shapes, it was recognized earlier that their surface energy was not that of surface tension for which only spherical shapes of closed membranes would occur. In 1973 the German physicist Wolfgang Helfrich (\cite{H}), based on previous studies about liquid crystallography, obtained the form of the energy which determines the geometry of these bilipid membranes
\begin{equation}\label{H}
\mathcal{H}_{a,c_o,b}[\Sigma]:=\int_\Sigma \left(a\left[H+c_o\right]^2+bK\right)d\Sigma\:,
\end{equation}
where $H$ and $K$ denote, respectively, the mean and Gaussian curvature of the mathematical surface $\Sigma$ representing the membrane. The Helfrich energy contains three parameters, which depend on the composition of the membrane itself: $a$ which is a bending rigidity modulus, $b$ the saddle-splay modulus and $c_o$ the spontaneous curvature (our definition for $c_o$ differs from the classical spontaneous curvature by the sign and a coefficient two). Assuming the membrane is homogeneous, these physical parameters are constant. As these constant parameters in the energy are slowly varied, the equilibrium states are expected to change in response. Bifurcations typically occur when one of these equilibria transitions between stability and instability. As this happens, new equilibria may suddenly appear or disappear.  

Biological membranes may have pores which allow materials to enter and leave cells. This introduces boundary components giving rise to equilibrium membranes with edges (\cite{MF,TOY,TOY03,TOY2,W}). Different studies of the Helfrich energy with an inelastic boundary line tension can be found in \cite{BR,C,CZT,Tu2,Z}, while those with elastic boundary can be found in \cite{AB,BMF,G17,PP1,PP2} and the references therein. The current investigation will concern a symmetry breaking bifurcation of a family of axially symmetric membranes having a fixed circular boundary. Previous studies of membrane bifurcation can be found in  \cite{BN,G17}. Since the boundary of the membrane will be kept fixed, the parameters $a$ and $b$ in the Helfrich energy \eqref{H} will not appear in the discussion below. To the contrary, the spontaneous curvature $c_o$ will arise in an essential way. We recall here that the spontaneous curvature $c_o$ originates primarily from the asymmetry between the two layers of the membrane (\cite{S}), although it may also arise from differences in the chemical properties of the fluids on both sides of the lipid bilayer (\cite{DSL}). Geometrically, this can be measured as the difference in area between the two layers, giving rise to the previous model for the Helfrich energy. The value of $c_o$ favors the local geometry of the membrane to be spherical, planar or hyperbolic (\cite{RL}).

At first, a rigorous approach to bifurcation in this case may seem prohibitively difficult due to the fact that the functional is quadratic in the surface curvatures, resulting in a fourth order Euler-Lagrange equation
\begin{equation}\label{EL} 
\Delta H+2\left(H+c_o\right)\left(H\left[H-c_o\right]-K\right)=0\,.
\end{equation} 
This would necessitate the linearized equation to have two boundary conditions in order to be well posed. Recently, the authors have found an argument which, in the case of an axially symmetric disc type surface, reduces the order of the Euler-Lagrange equation \eqref{EL} to two (\cite{PP2}). This lower order functional is equivalent to the problem of finding equilibria for a linear combination of the surface area and the gravitational potential energy of the surface when it is regarded as lying in the three dimensional hyperbolic space (\cite{PP2}).

In this paper we will study a symmetry breaking bifurcation of a specific one parameter family of axially symmetric membrane shapes. To produce this bifurcation, Crandall-Rabinowitz theory is applied to specific solutions of the second order problem mentioned above. Our method is closely related to that used in \cite{P} to produce a symmetry breaking bifurcation of a one parameter family of constant mean curvature surfaces. In fact, the bifurcation of the family of nodoids produced in \cite{P} could be considered a symmetry breaking bifurcation of annular membranes since constant mean curvature surfaces trivially solve the equation \eqref{EL} when the value of the mean curvature is taken to be $H\equiv -c_o$. We note that in the case considered in this paper, the mean curvature of the surfaces is not constant.

We begin by proving in Theorem \ref{Sigma0} that for a given circle, which will represent the surface boundary, there exists an axially symmetric disc type surface $\Sigma_0$ (see the figure on the front page) spanning this circle and satisfying the second order reduction of the Euler-Lagrange equation \eqref{EL} obtained in \cite{PP2}, for some value of the spontaneous curvature $c_o$. Moreover, we show that this surface meets the boundary circle tangentially. In a second step, Theorem \ref{A}, we show, by means of the Implicit Function Theorem, the existence of a one parameter family of axially symmetric disc type surfaces $\Sigma_t$, $t\in(-\epsilon,\epsilon)$ for sufficiently small $\epsilon>0$, with the same boundary circle as $\Sigma_0$ but with varying contact angle. This family is illustrated in Figure \ref{family}. Finally, in Theorem \ref{bifurcation}, from an argument involving simple eigenvalues (\cite{CR}) we deduce the existence of a non-axially symmetric family of surfaces bifurcating from $\Sigma_t$ at, precisely, $\Sigma_0$. A first order linear approximation of a branch of this family is shown in Figure \ref{NewFamily}.

\section{Preliminaries}

The membrane will be modeled as a compact, connected and oriented surface $\Sigma$ which we regard as being smoothly embedded in the three dimensional Euclidean space 
$$X:\Sigma \rightarrow {\bf R}^3\:.$$
The canonical coordinates of ${\bf R}^3$ will be represented by $(x,y,z)$, while $E_i$, $i=1,2,3$ will be the constant unit vector fields in the direction of the coordinate axes. We denote the unit normal field by $\nu:\Sigma \rightarrow {\bf S}^2$, defined by the property that it points out of any convex domain. In this paper, the surface $\Sigma$ will always be a topological disc and the embedding is assumed to be smooth up to the boundary.

In a previous paper (\cite{PP2}), we showed that any axially symmetric disc type surface which is critical for a functional $\mathcal{H}_{a,c_o,b}$ is necessarily critical for a lower order functional. We repeat some of these results here for the sake of completeness.

\begin{theorem}\label{prop1}\mbox{\cite[Theorem~4.1]{PP2}} If $\Sigma$ is a smooth disc type axially symmetric surface satisfying the equation \eqref{EL}, then either the surface has constant mean curvature or an equation
\begin{equation}\label{nonCMC}
H+c_o=\frac{\nu_3}{A-z}\,,
\end{equation}
holds for a suitable constant $A$. (By making the linear change of the vertical coordinate, we can assume $A=0$.)
\end{theorem}

The condition \eqref{nonCMC} is the Euler-Lagrange equation for the functional (see Theorem 4.2 and Remark 4.1 of \cite{PP2})
\begin{equation}\label{G}
\mathcal{G}[\Sigma]:=\widetilde{\mathcal{A}}[\Sigma]-2c_o\int_{\widetilde{V}}\lvert z\rvert\,d\widetilde{V}=\widetilde{\mathcal{A}}[\Sigma]-2c_o\mathcal{U}[\Sigma]\,,
\end{equation}
where $\widetilde{A}$ denotes the area of $\Sigma$, regarded as a surface in ${\bf H}^3$, $\widetilde{V}$ denotes the hyperbolic volume enclosed by $\Sigma$ and $\mathcal{U}$ is the gravitational potential energy of $\Sigma$, considered as a surface in ${\bf H}^3$. In this variational problem we consider only those variations of the surface vanishing on the boundary. In terms of the metric on ${\bf R}^3$,
$$\mathcal{G}[\Sigma]=\int_\Sigma \frac{1}{z^2}\,d\Sigma+2c_o\int_\Sigma\frac{\nu_3}{z}\,d\Sigma\,,$$
after applying the Divergence Theorem. 

The converse of Proposition \ref{prop1} holds, whether the surface is axially symmetric or not and regardless of the topological type of the surface.  For convenience, we will translate the vertical coordinate so that the constant $A$ in \eqref{nonCMC} is zero. 

\begin{prop}\mbox{\cite[Proposition~4.1]{PP2}} Assume that the equation
\begin{equation}\label{H+c}
H+c_o=-\frac{\nu_3}{z}\,,
\end{equation}
holds on $\Sigma$. Then, \eqref{EL} also holds. 
\end{prop}

In the axially symmetric case, surfaces satisfying the relation \eqref{H+c} can be generated by solving the system for the arc length parameterized generating curve $\gamma(\varsigma)=\left(r(\varsigma),z(\varsigma)\right)$ given by
\begin{eqnarray}\label{odesys}
r_\varsigma(\varsigma)&=&\cos\varphi(\varsigma)\,,\nonumber\\
z_\varsigma(\varsigma)&=&\sin\varphi(\varsigma)\,,\\
\varphi_\varsigma(\varsigma)&=&-2\frac{\cos\varphi(\varsigma)}{z(\varsigma)}-\frac{\sin\varphi(\varsigma)}{r(\varsigma)}+2c_o\,,\nonumber
\end{eqnarray}
where we are denoting by $\left(\,\right)_\varsigma$ the derivative with respect to the arc length parameter $\varsigma$. Here, the function $\varphi(\varsigma)$ represents the angle between the positive part of the $r$-axis and the tangent vector to $\gamma(\varsigma)$.

\begin{remark} As mentioned in Remark 5.1 of \cite{PP2}, equation \eqref{H+c} carries both the positive and negative signs in front of $c_o$ in above system. However, up to the transformation $z\mapsto -z$, this sign can be fixed to be positive. Similarly, we may assume $c_o\geq 0$.
\end{remark}

In this paper, we will apply this system to generate disc type surfaces with $\varsigma=0$ corresponding to the boundary of the surface and $\varsigma =\ell$ corresponding to the cut with the axis of rotation. Since we are looking for disc type surfaces, the generating curve $\gamma(\varsigma)$ must cut the axis of rotation and so $r(\ell)=0$ holds. Of course, the value of $\ell$ will depend on the individual surface and will determine the length of the generating curve. Moreover, for regularity, we also need to request the cut with the rotation axis to be perpendicular. This completely describes the tangent vector at that point, which is equivalent to fixing the value $\varphi(\ell)=\pi$. 

From Remark 5.2 of \cite{PP2} we conclude that, when the spontaneous curvature $c_o$ is zero, the only solutions of \eqref{odesys} with above conditions at $\varsigma=\ell$ are parts of circles and so the corresponding surfaces are spherical caps. Therefore, we assume from now on that $c_o$ is positive. For fixed spontaneous curvature $c_o>0$ we will have a one parameter family of axially symmetric disc type surfaces. The parameter of this family corresponds with the value of the height at $\varsigma=\ell$, i.e., $z(\ell)=z_o$.

The system \eqref{odesys} is singular at $r=0$ and, hence, the existence of solution for the conditions at $\varsigma=\ell$ described above is not guaranteed by the standard theory of ordinary differential equations. Nevertheless, in the following result, we will prove the existence of classical solutions.

\begin{prop}\label{ex} Let $c_o>0$ and $z_o\neq 0$. Then, the system of first order differential equations \eqref{odesys} with $r(\ell)=0$, $z(\ell)=z_o$ and $\varphi(\ell)=\pi$ has analytic solutions. Moreover, these solutions depend continuously on the parameters $c_o$ and $z_o$.
\end{prop}
{\it Proof.\:} For this proof, we will follow the techniques of \cite{ODEPaper}. Near $r=0$, we rewrite the system \eqref{odesys} by describing the generating curve $\gamma$ as a graph $z\equiv z(r)$,
\begin{equation}
\label{nonpar}\left(\frac{r z_r}{2\sqrt{1+z_r^2}}\right)_r=-\frac{r}{z\sqrt{1+z_r^2}}-c_o r\,\end{equation}
In this notation the initial conditions reduce to $z(0)=z_o$ and $z_r(0)=0$.

The local existence and uniqueness of solutions of this initial value problem can be deduced from the fixed points of a suitable operator. Indeed, define the operator $\mathcal{T}$ acting on $z$ as
$$\left(\mathcal{T}\,z\right)[r]:=z_o+\int_0^r \Psi^{-1}\left(\int_0^s-2\frac{t}{s}\left[\frac{1}{z\sqrt{1+z_t^2}}+c_o\right]dt\right)ds\,,$$
where $\Psi(u):=u/\sqrt{1+u^2}$. It is then clear that fixed points of $\mathcal{T}$ give rise to solutions of above initial value problem.

The inverse function of $\Psi$ is given by $\Psi^{-1}(t):=t/\sqrt{1-t^2}$. Differentiating we obtain that 
$$\left[\Psi^{-1}\right]'(t)=\left(1-t^2\right)^{-3/2}.$$  
Since $z_o\neq 0$ and $s$ can be chosen in the inner integral, the function $\Psi$ is globally Lipschitz continuous. We can then check that the operator $\mathcal{T}$ is a contraction in the space $\mathcal{C}^1([0,\epsilon])$ for $\epsilon>0$ sufficiently small. Therefore, by the Contraction Mapping Theorem, for fixed $c_o$ and $z_o$ there exists a unique fixed point of $\mathcal{T}$, and so a unique local solution of the initial value problem. The continuous dependence on the parameters $c_o$ and $z_o$ of these local solutions follows from the continuous dependence of the fixed points of the operator $\mathcal{T}$.

Far from $r=0$ the standard theory of ordinary differential equations can be applied to extend the solution and the analyticity is a consequence of the elliptic regularity. {\bf q.e.d.}
\\

\begin{figure}[h!]
\begin{center}
\includegraphics[width=0.6\linewidth]{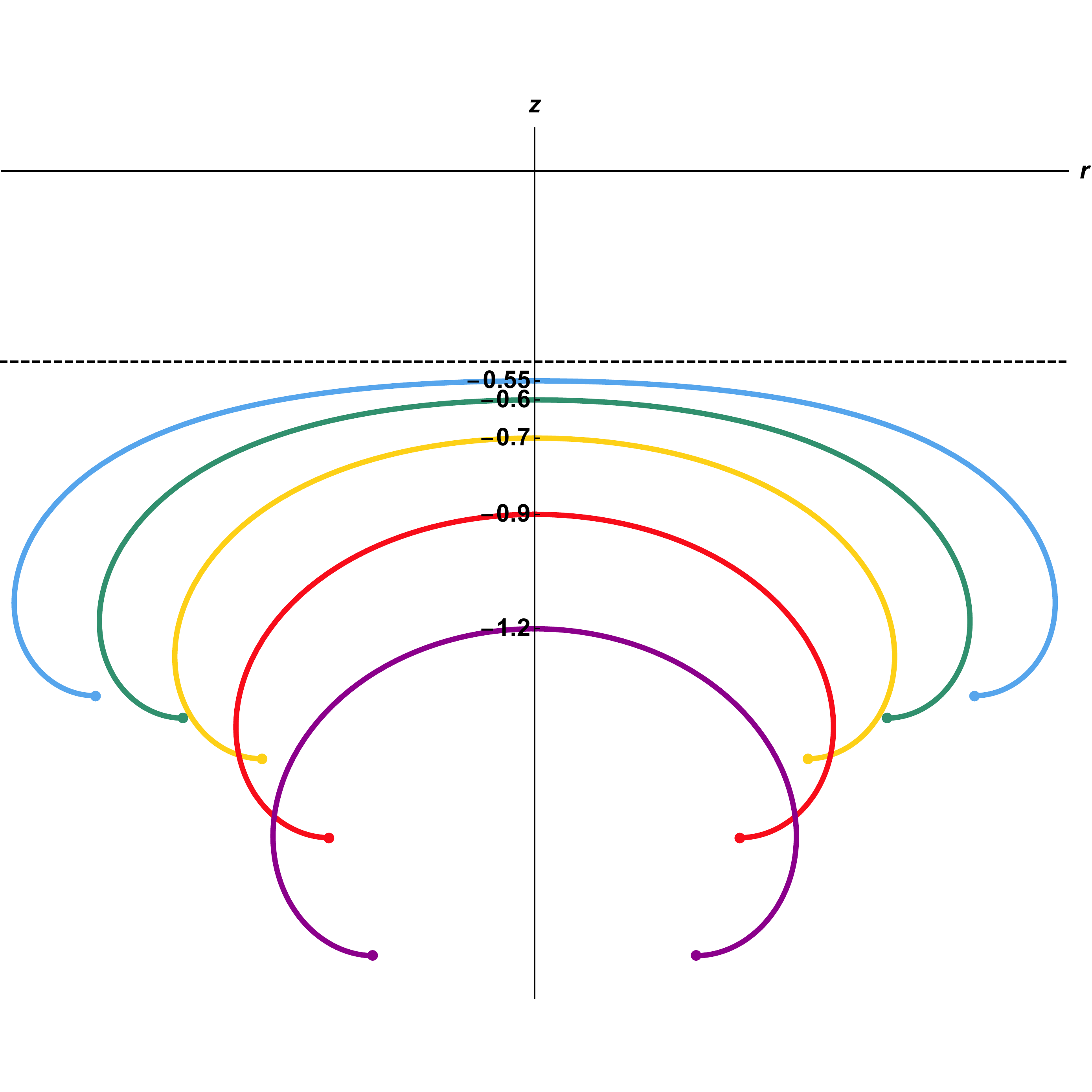}
\end{center} 
\caption{Several generating curves $\gamma$ for the corresponding parameters $z_o<-1/c_o$. The dashed line represents the limit line $z=-1/c_o$. Here, $c_o=2$.}
\label{Profiles}
\end{figure}

A particular example, which will be central to our discussion, is the surface $\Sigma_0$ shown on the front page. The tangent plane along the boundary circle of this surface is horizontal and $\varphi(0)=0$. Surfaces of this type appear provided that $z_o<-1/c_o<0$. Indeed, we will show in what follows that for the parameter $z_o$ satisfying such a relation all the generating curves $\gamma(\varsigma)$ have the same shape (see Figure \ref{Profiles}), and so all the corresponding axially symmetric disc type surfaces are of the same type of $\Sigma_0$. 

\begin{theorem}\label{shape} Let $c_o>0$ and $z_o<0$ be constants satisfying $z_o<-1/c_o$. Then, there exists a curve $\gamma(\varsigma)=\left(r(\varsigma),z(\varsigma)\right)$, $\varsigma\in[0,\ell]$, which is a solution of the system \eqref{odesys} with $r(\ell)=0$, $z(\ell)=z_o$ and $\varphi(\ell)=\pi$, such that, when rotated about the vertical axis, the resulting surface $\Sigma_0$ is a convex topological disc for which $z_o$ is its maximum height. Moreover, $\Sigma_0$ is bounded by a circle along which the tangent plane is horizontal.
\end{theorem}
{\it Proof.\:} Let $\gamma(\varsigma)$ be a solution of the system \eqref{odesys} with $r(\ell)=0$, $z(\ell)=z_o$ and $\varphi(\ell)=\pi$. A short time solution of this type is guaranteed by Proposition \ref{ex}.

First, we will see that the angle $\varphi$ is increasing near $\varsigma=\ell$ and we will employ this to geometrically describe $\gamma$ near the cut with the vertical axis. Using that $\varphi(\ell)=\pi$ and L'Hopital's rule in the last equation of \eqref{odesys}, we obtain that
$$\varphi_\varsigma(\ell)=-2\frac{\cos\varphi(\ell)}{z(\ell)}-\lim_{\varsigma\to\ell}\frac{\sin\varphi(\varsigma)}{r(\varsigma)}+2c_o=\frac{2}{z_o}-\varphi_\varsigma(\ell)+2c_o\,,$$
and so $\varphi_\varsigma(\ell)=1/z_o+c_o>0$, due to the relation between $c_o>0$ and $z_o<0$. Consequently, $\varphi$ increases for $\varsigma\in\left(\ell-\epsilon,\ell\right)$ with $\epsilon>0$ sufficiently small. In particular, it follows that $\varphi(\varsigma)<\pi$, but close to it, when $\varsigma\in(\ell-\epsilon,\ell)$ and, hence, the first two equations of \eqref{odesys} imply that $z_\varsigma=\sin\varphi>0$ and $r_\varsigma=\cos\varphi<0$, that is, $z$ increases while $r$ decreases. Since we are going backwards this means that $z(\varsigma)\leq z_o$ and $r(\varsigma)\geq 0$ for every $\varsigma\in\left(\ell-\epsilon,\ell\right]$. Observe that the curve $\gamma(\varsigma)$ can be continued to all values of $\varsigma<\ell$ as long as $z<0$ and $r>0$ hold.

In a second step, we will show that a vertical tangent vector of the system \eqref{odesys} occurs. To this end, we will employ the nonparametric representation of the curve $\gamma$ as a graph $z\equiv z(r)$, which is given by \eqref{nonpar}. Comparing \eqref{nonpar} with the system \eqref{odesys}, we obtain
$$-\left(r\sin\varphi\right)_r=-\sin\varphi-r\varphi_\varsigma=-2\left(-\frac{\cos\varphi}{z}+c_o\right)r\leq -2\left(\frac{1}{z_o}+c_o\right)r\,.$$
Integrating the previous inequality with respect to $r$, from $0$ to $r$, and using the initial conditions of \eqref{nonpar}, yields
$$\sin\varphi\geq\left(\frac{1}{z_o}+c_o\right)r\,.$$
Since $1/z_o+c_o>0$ holds by assumption, one sees that a vertical tangent must be reached before $r$ reaches the value $(1/z_o+c_o)^{-1}$. We will denote by $\widehat{\varsigma}$ the value of the arc length parameter where the vertical tangent is reached and by $\overline{\gamma}$ the part of the curve $\gamma$ which is the image of the interval $(\widehat{\varsigma},\ell]$.

Next, we will show that $\overline{\gamma}$ is convex. This will follow from the monotonicity of $\varphi$ in that part of the curve. We know that $\varphi$ is increasing for $\varsigma\in(\ell-\epsilon,\ell]$. If $\varphi$ were to stop increasing in $\overline{\gamma}$, there would exist a value $\widetilde{\varsigma}\in(\widehat{\varsigma},\ell)$ such that $\varphi_\varsigma(\widetilde{\varsigma})=0$ and $\varphi_{\varsigma\varsigma}(\widetilde{\varsigma})\geq 0$. However, assume that $\varphi_\varsigma(\widetilde{\varsigma})=0$ holds, then differentiating $\varphi_\varsigma$ from \eqref{odesys} and evaluating at $\widetilde{\varsigma}$ we get
$$\varphi_{\varsigma\varsigma}(\widetilde{\varsigma})=\sin\varphi(\widetilde{\varsigma})\cos\varphi(\widetilde{\varsigma})\left(\frac{2}{z^2(\widetilde{\varsigma})}+\frac{1}{r^2(\widetilde{\varsigma})}\right),$$
which is negative since $\varphi\in(\pi/2,\pi)$. This is a contradiction, so $\varphi$ is always increasing in $\overline{\gamma}$. Therefore, the curvature of $\gamma$ satisfies $\kappa=-\varphi_\varsigma<0$. Up to reversing the orientation, this shows that $\overline{\gamma}$ is convex.

At the point where the vertical tangent is reached, $\nu_3(\widehat{\varsigma})=-r_\varsigma(\widehat{\varsigma})=0$. Since $\nu_3$ satisfies the second order differential equation \eqref{Pnu3} (the validity of this equation will be shown below), it has a sign change and below the points where $\nu_3=0$, $\nu_3<0$ holds. This means that $r(\varsigma)$ has a local maximum at $\widehat{\varsigma}$. In some interval $(\widehat{\varsigma}-\epsilon,\widehat{\varsigma})$ for $\epsilon>0$ sufficiently small, the angle $\varphi$ lies in the interval $(0,\pi/2)$. We denote by $\underline{\gamma}$ the image of this interval.

We now show that $r$ is bounded below by a positive constant for any $\varsigma<\ell$. We multiply the third equation in \eqref{odesys} by $\cos\varphi$ to obtain
$$\left(\sin \varphi\right)_\varsigma =\varphi_\varsigma \cos \varphi =-2\frac{\cos^2\varphi}{z}-\frac{\cos \varphi\sin\varphi}{r}+2c_o\cos\varphi\:.$$
Multiplying by $r$ and integrating from $\varsigma$ to $\ell$, gives
$$\int_\varsigma ^\ell (\sin \varphi)_\varsigma r\:d\varsigma = \int_\varsigma ^\ell \frac{-2\cos^2\varphi}{z} r\:d\varsigma-\int_\varsigma ^\ell\cos \varphi\sin\varphi\:d\varsigma+ 2c_o\int_\varsigma ^\ell\cos\varphi r\,d\varsigma\:.$$
Integrating the first term by parts and using $\sin\varphi=z_\varsigma$ and $\cos\varphi=r_\varsigma$, we arrive at
$$c_or^2-r\sin \varphi - \int_\varsigma ^\ell \frac{-2\cos^2\varphi}{z} r\:d\varsigma=0\:.$$
So
$$r=\frac{1}{2c_o}\left(\sin\varphi \pm \left[\sin^2\varphi+4c_o\int_\varsigma ^\ell \frac{-2\cos^2\varphi}{z}r\:d\varsigma\right]^{1/2}\right).$$ 
As long as $z<0$ holds, the integral in the last expression is positive. It follows that, since the discriminant is  positive on any set  $\varsigma <\ell$,  $r$ is bounded below by a positive constant on any set where $\varsigma <\ell$ and only the root with the positive sign can occur.  

We next see that $\underline{\gamma}$ is convex. For this purpose we will use another consequence of above formula for $r$, which gives that
\begin{equation*}
c_o-\frac{\sin \varphi}{r} \ge 0\,,
\end{equation*}
holds on $\gamma$. This implies that  any arc of $\gamma$ on which $\cos \varphi>0$ holds is convex. This follows
since the curvature of $\gamma$ is
$$\kappa =-\varphi_\varsigma =\frac{2\cos\varphi}{z}+\frac{\sin \varphi}{r}-2c_o< 0\:.$$

We finish showing that the angle $\varphi=0$ is attained. If the arc ${\underline \gamma}$ contains a point where its tangent is horizontal, then we are done so assume no such point exists. The intersection of any tangent line to ${\underline \gamma}$ with  the vertical axis gives a lower bound for $z$ among points in ${\underline \gamma}$. Let ${\underline z}={\rm inf}_{{\underline \gamma}}z $. The arc ${\underline \gamma}$ can be represented as a graph $r=r(z)$ satisfying the equation
$$\left(\frac{r_z}{\sqrt{1+r_z^2}}\right)_z=\frac{2r_z}{z\sqrt{1+r_z^2}}+\frac{1}{r\sqrt{1+r_z^2}}-2c_o\:.$$
Since there are no horizontal tangents, we get that $r_z$ is bounded and so the solution can be extended to an open neighborhood of ${\underline z}$ which is a contradiction. {\bf q.e.d.}
 \\
 
We now prove the existence of a surface of type $\Sigma_0$ for any circle, which serves as the boundary curve.

\begin{theorem}\label{Sigma0} For every circle $C$ there exists an axially symmetric convex embedded disc type surface $\Sigma_0$ satisfying
$$H+c_o=-\frac{\nu_3}{z}\,,$$
for some value of $c_o$ and having a horizontal tangent plane along its boundary $C$.
\end{theorem}
{\it Proof.\:} After the transformation $z\mapsto -z$, if needed, we assume that the circle $C$ is located in the semi-space $z<0$, and that $c_o>0$ holds.

The generating curve of an axially symmetric disc type surface in this setting is a solution of the system \eqref{odesys} for $r(\ell)=0$, $z(\ell)=z_o$ and $\varphi(\ell)=\pi$. The existence of such a solution is guaranteed by Proposition \ref{ex}. We just consider the family of solutions arising for $z_o<-1/c_o$. For these solutions, as proved in Theorem \ref{shape}, we can restrict the length of the curves $\ell$ such that $\varphi(0)=0$.

Observe that when $c_o\to\infty$ and $z_o\to 0$, the solution approaches the straight line $z=0$, while when $c_o\to\infty$ and any $z_o<0$, the solution tends to a semicircle. From the continuous dependence of solutions of \eqref{odesys} on the parameters $c_o$ and $z_o$ (see Proposition \ref{ex}), by varying these parameters in their domains $c_o>0$ and $z_o<-1/c_o$ we cover all possible radii and heights, so we prove the existence of solution for a given circle. Moreover, by the restriction on the value $\ell$, this solution meets the boundary circle at an angle $\varphi(0)=0$. {\bf q.e.d.}

\section{Linearization and Eigenvalue Problem}

Our first objective is to produce a one parameter family of axially symmetric surfaces $\{\Sigma_t\}$, $t\in(-\epsilon,\epsilon)$ with $\epsilon>0$ sufficiently small, all satisfying an equation \eqref{H+c} with $c_o\equiv c_o(t)$ and all having the same boundary as $\Sigma_0$. This will be achieved by using the Implicit Function Theorem (\cite{Heida}), so we will first need to discuss the linearization of the equation \eqref{H+c} on $\Sigma_0$. It will be convenient to define, on any surface on which $z<0$, a function $\xi$ by
\begin{equation}\label{xi}
\xi:=H+\frac{\nu_3}{z}\:.
\end{equation}
Clearly, the equation \eqref{H+c} characterizing our surfaces can be expressed as $\xi=-c_o$.
  
We express a variation of the immersion $X_0:\Sigma_0\rightarrow{\bf R}^3$ as $\delta X_0=\psi\nu +(\delta X_0)^T$ with $(\delta X_0)^T$ tangent to $\Sigma_0$ and we then define an operator $P$ using the pointwise variation of $\xi$ by
\begin{equation}\label{dxi}
\delta\xi=\frac{1}{2}P[\psi]+\nabla \xi\cdot (\delta X_0)^T=\frac{1}{2}P[\psi]\:,
\end{equation}
since $\xi$ is constant on $\Sigma_0$. In order to specify the operator $P$, recall that the variation of the mean curvature for the normal variation $\delta X_0=\psi\nu$ is (for details, see Appendix A of \cite{PP2})
$$\delta H=\frac{1}{2}L[\psi]=\frac{1}{2}\left(\Delta \psi+\left[4H^2-2K\right]\psi\right).$$ 
This gives, from \eqref{xi}, that
$$P[\psi]=L[\psi]+2\frac{z\delta \nu_3-\nu_3\delta z}{z^2}=L[\psi]-2\frac{\nabla \psi\cdot \nabla z}{z}-2\frac{\nu^2_3}{z^2}\psi\:.$$
Therefore,
\begin{equation*}
P[\psi]= z^2\left(\nabla\cdot\left[z^{-2}\nabla \psi\right]+U\psi \right),
\end{equation*}
where $U:=z^{-2}\left(4H^2-2K-2[\nu_3/z]^2\right)$. For a critical surface, the linearization of equation \eqref{H+c} is then $P[\psi]=0$. \emph{This assumes that the value of $c_o$ is held fixed throughout the variation}. Since any horizontal translation is a symmetry to the variational problem \eqref{G}, it follows that $P[\nu_i]=0$, $i=1,2$. Note that $\nu_i=0$ on $\partial\Sigma_0$.

More generally, we can consider the eigenvalue problem
\begin{equation}\label{P}
P[\psi]+\lambda z^2 \psi=0\,,\quad\quad\quad \psi\lvert_{\partial \Sigma_0}\equiv 0\:.
\end{equation}
If the surface $\Sigma_0$ is regarded as the image of $X_0(\varsigma,\theta)=\left(r(\varsigma)e^{i\theta}, z(\varsigma)\right)$, then \eqref{P} becomes
\begin{equation}\label{sep1}
\frac{1}{r}\left[\frac{r\,\psi_\varsigma}{z^2}\right]_\varsigma+\frac{\psi_{\theta \theta}}{r^2z^2} +U(\varsigma) \psi+\lambda  \psi=0\,, \quad\quad\quad u(0)=0\,.\end{equation}
In particular, if we separate variables and write $\psi(\varsigma,\theta)=u(\varsigma)\cos[m\theta]$ (or, equivalently, $\psi(\varsigma,\theta)=u(\varsigma)\sin [m\theta]$), then $u$ must satisfy
\begin{equation}\label{sep}
\frac{1}{r}\left[\frac{r\,u_\varsigma}{z^2}\right]_\varsigma-\frac{m^2u}{r^2z^2} +U(\varsigma) u+\lambda  u=0\,, \quad\quad\quad u(0)=0\,.
\end{equation}
{\it Note that if $m\ge 1$, then $u(\ell)=0$ is necessary for regularity at $\varsigma=\ell$.}

\begin{lemma} There exists a $\mathcal{C}^2$ axially symmetric solution of the differential equation $P[\psi(\varsigma)]=0$ on $\Sigma_0$. The dimension of the space of axially symmetric solutions of this equation is one.
\end{lemma}
{\it Proof.\:} From \eqref{sep1} with $\lambda=0$, the equation $P[\psi(\varsigma)]=0$ becomes
\begin{equation*}
\frac{1}{r}\left[\frac{r\,\psi_\varsigma}{z^2}\right]_\varsigma+U(\varsigma)\psi=0\:,
\end{equation*}
which is singular at $\varsigma=\ell$ since $r(\ell)=0$.

We first produce an axially symmetric solution on an interval $[\widetilde{\varsigma}, \ell]$. Let $\Omega$ be an axially symmetric domain containing the
top of the surface for which the first eigenvalue of the problem $P[f]+\lambda f=0$, $f\lvert_{\partial \Omega}=0$ is positive. By the Fredholm alternative, $\Omega$ has a solution of $P[g]=0$, $g\lvert_{\partial \Omega}=1$. Since the coefficients of $P$ only depend on $\varsigma$, this equation can be differentiated with respect to $\theta$, to get $P[g_\theta]=0$ and $g_\theta\lvert_{\partial \Omega}=0$. So, since $g_\theta$ is an eigenfunction 
for $0$ in $\Omega$ and the first Dirichlet eigenvalue is positive, it follows that $g_\theta\equiv 0$ and $g=g(\varsigma)$ is an axially symmetric solution. 

The solution of the system of ordinary differential equations \eqref{odesys} is smooth and real analytic on $(0,\ell)$ and so are the coefficients of the operator $P$. By standard theory of ordinary differential equations, the function $g$ can be smoothly extended to all of $\Sigma_0$. {\bf q.e.d.}

\begin{lemma}\label{nas} If $\psi(\varsigma)\not\equiv 0$ is an axially symmetric solution of $P[\psi]=0$, then $\psi\lvert_{\partial \Sigma_0}\neq 0$.
\end{lemma}
{\it Proof.\:} Assume to the contrary of the statement of the lemma that $\psi(0)=0$. Since $P[\nu_1]=0$ and $\nu_1$ has a sign change in $\Sigma_0$, $\lambda=0$ is not the first eigenvalue of the problem \eqref{P} and so $\psi(\widetilde{\varsigma})=0$ for some $\widetilde{\varsigma}\in (0,\ell)$.

Consider the variation of the embedding $X_0:\Sigma_0\rightarrow {\bf R}^3$, given by $X_\rho:= X_0+\rho E_3$. The variation field for this variation has normal component $\partial_\rho (X_\rho)_{\rho=0}\cdot \nu =\nu_3$. For this variation 
$$\xi_\rho =H+\frac{\nu_3}{z+\rho}\:,$$
so it follows from \eqref{dxi} that
\begin{equation}\label{Pnu3}
P[\nu_3]+2\frac{\nu_3}{z^2}=0\:,
\end{equation}
holds. If ${\widehat\varsigma}$ denotes the unique value in $[0,\ell]$ with $\nu_3({\widehat \varsigma})=0$, then it follows from the Sturm Comparison Theorem that ${\widehat \varsigma}<\widetilde{\varsigma}$ holds. There may possibly be other zeroes of $\psi$ but we assume that $\widetilde{\varsigma}$ is the largest value of $\varsigma<\ell$ with $\psi(\widetilde{\varsigma})=0$ and we let $\Omega$ be the domain of $\Sigma_0$ bounded by the circle $\varsigma=\widetilde{\varsigma}$. Note $\nu_3>0$ holds on $\Omega$.

We assume $\psi\geq 0$ in $\Omega$ (if not, we replace $\psi$ with $-\psi$ and rename). The operator $P$ is self adjoint with respect to the measure $z^{-2}\:d\Sigma$ or equivalently $z^{-2}P$ is self adjoint. Using this,  we get
$$0>-2\int_\Omega \psi\frac{\nu_3}{z^4}\:d\Sigma=\int_\Omega \psi P[\nu_3] z^{-2}\:d\Sigma=\int_\Omega  \left(\psi P[\nu_3] -\nu_3P[\psi]\right)z^{-2}\:d\Sigma\,,$$
where we have used \eqref{Pnu3} and that $P[\psi]$ is assumed to be zero. Now, from the expression of $P$ and after integrating by parts, we obtain that the integral on the right-hand side above equals
$$\oint_{\partial \Omega}\left(\psi z^{-2}\nabla\nu_3-\nu_3z^{-2} \nabla \psi \right)\cdot n\:ds= -2\pi\, r(\widetilde{\varsigma})\left(\nu_3z^{-2} \nabla \psi \cdot n\right)_{\partial \Omega}>0\:,$$
where we have used that $\psi(\widetilde{\varsigma})=0$ and with the last inequality following since $\psi\not\equiv 0$ and $\nabla\psi\cdot n<0$. This gives a contradiction. {\bf q.e.d.}

\begin{lemma}\label{m} For $m>1$ there is no non trivial solution of \eqref{sep} vanishing on $\partial \Sigma_0$.
\end{lemma}
{\it Proof.\:} As noted above, any eigenfunctions of the problem \eqref{P} can be found by separation of variables. If $\psi(\varsigma,\theta)=u(\varsigma)\cos [m\theta]$ or $\psi(\varsigma,\theta)=u(\varsigma)\sin [m\theta]$ is an eigenfunction with $m\geq 1$, then it is clear that $u(\ell)=0$, otherwise $\lim_{\varsigma \rightarrow \ell}\psi$ does not exist, (recall that $\ell$ is the value of $\varsigma$ at the top of $\Sigma_0$).

Since $z_\varsigma$ does not vanish in the interior of $\Sigma_0$, it follows that $u=\zeta z_\varsigma$ for a well defined function $\zeta$. For convenience, we set $w:=z_\varsigma$ and note that $w$ solves \eqref{sep} with $m=1$, which is a direct consequence of $P[\nu_1]=0$. If $u$ is as above, we get
\begin{eqnarray*}
0&=&\int_0^\ell u\left( \frac{1}{r}\left[\frac{r u_\varsigma}{z^2}\right]_\varsigma -\frac{m^2u}{r^2z^2}+Uu\right)r\,d\varsigma\\
&=&\int_0^\ell u\left(\left[\frac{r(\zeta w)_\varsigma}{z^2}\right]_\varsigma -\frac{m^2 u}{rz^2}+rU u\right)d\varsigma\\ 
&=&\int_0^\ell u\left(\left[\frac{r(\zeta w_\varsigma+w\zeta_\varsigma)}{z^2}\right]_\varsigma -\frac{m^2u}{rz^2}+rUu\right)d\varsigma\\ 
&=&\int_0^\ell u\left(\zeta\left[\frac{r w_\varsigma}{z^2}\right]_\varsigma+2\frac{r\zeta_\varsigma w_\varsigma}{z^2}+w\left[\frac{r \zeta_\varsigma}{z^2}\right]_\varsigma -\frac{m^2u}{rz^2}+rUu\right)d\varsigma\\ 
&=&\int_0^\ell u\left([1-m^2]\frac{u}{rz^2}+2\frac{r\zeta_\varsigma w_\varsigma}{z^2}+w\left[\frac{r \zeta_\varsigma}{z^2}\right]_\varsigma\right)d\varsigma\,,
\end{eqnarray*}
where in the last equality we have used that $w$ is a solution of \eqref{sep} for $m=1$. Integrating the last term by parts, we get
$$0=\left(\frac{r\zeta_\varsigma w}{z^2}\right)_0^\ell+\int_0^\ell \left( [1-m^2]\frac{u^2}{rz^2}-\frac{\zeta_\varsigma^2w^2r}{z^2}\right)d\varsigma<0\:.$$
Both the functions $w$ and $u$ satisfy the second order equation \eqref{sep} and vanish at $\varsigma=0$, $\ell$. The regularity of $\zeta$ at these points can be deduced from the local behavior of the zeros of elliptic equations, so the first term in the last equation vanishes and a contradiction is reached. {\bf q.e.d.}
\\
 
Let $q:=X_0\cdot \nu$ be the support function of $\Sigma_0$. If the embedding $X_0$ is rescaled $X_0\mapsto (1+\rho)X_0$, then $\xi:=H+\nu_3/z  \mapsto (1+\rho)^{-1}\xi=-(1+\rho)^{-1}c_o$. It follows that
$$\frac{1}{2}P[q]=\partial_\rho \left(\frac{\xi}{1+\rho}\right)_{\rho=0}=c_o\:.$$

By Lemmas \ref{nas} and \ref{m}, the only axially symmetric solution of $P[\psi]=0$, $\psi\lvert_{\partial \Sigma_0}=0$, is $\psi\equiv 0$. It is clear that $q$ does not vanish on $\partial \Sigma_0$ since the normal to $\Sigma_0$ along the boundary is the vector $(0,0,-1)$. So 
\begin{equation}
\label{h}h:=c_o^{-1}\left(q(0)\psi-\psi(0)q\right),
\end{equation}
satisfies  
\begin{equation*}
P[h]=-2, \quad h\lvert_{\partial \Sigma_0}=0\:.
\end{equation*}

We are now in a position to exert the existence of a one parameter family of axially symmetric membranes containing $\Sigma_0$ which all share the same boundary. (Cf. \cite{K}, Theorem 1.2.)

\begin{theorem}\label{A}Let $\mathcal{A}_o^{k+\alpha}$ be the Banach space of axially symmetric functions in $\mathcal{C}_o^{k+\alpha}(\Sigma_0)$.
There exists a neighborhood $\widetilde{I}$ of $0\in {\bf R}$ and a $\mathcal{C}^2$ curve $\phi:{\widetilde I}\rightarrow \mathcal{A}_o^{k+\alpha}$ with 
$\phi(0)=0$, $\partial_t\phi(0)=h$ and such that $X(t)=X_0+\phi(t)\nu$ defines a one parameter family of axially symmetric surfaces $\{\Sigma_t\}$, $t\in \widetilde{I}$, satisfying the equation
$$\xi\equiv \xi(t)=H+\frac{\nu_3}{z}=c_o+t\:.$$ 

\end{theorem}
{\it Proof.\:} For a sufficiently small neighborhood of $I\subset \bf{R}$  containing $0$ and a sufficiently small neighborhood $U$ of $0$ in $\mathcal{A}_o^{k+\alpha}$, define 
$$\Xi: I \times U \rightarrow \mathcal{A}^{k-2+\alpha}, \quad \Xi(t,f):=\xi\left(X_0+[th+f]\nu\right)+(c_o+t)\:.$$
Then $\Xi(0,0)=0$ and the partial derivative of this map with respect to $f$ at $t=0$ and $f=0$ applied to $\zeta \in \mathcal{A}_o^{k+\alpha} $ is
$$\left(D_f\Xi \right)_{(0,0)}[{\zeta}]:=\partial_\epsilon \left(\xi(X_0+[0+\epsilon \zeta]\nu\right)+c_o)=P[\zeta]\:. $$
Although it is cumbersome to compute the derivatives explicitly, it can be shown that the map $\Xi$ is, in fact, of class $\mathcal{C}^2$.

Under the assumption given above, this derivative is non singular, so by the Implicit Function Theorem (\cite{Heida}), there exists a $\mathcal{C}^2$ curve $f(t)$ from an interval $0\in \widetilde{I}\subset I$ into $\mathcal{A}^{k+\alpha}_o$ such that $\Xi(t, th+f(t))\equiv 0$. For $t\approx 0$, the surfaces $X_0+(th+f(t))\nu$ comprise a one parameter family of axially symmetric surfaces sharing the same boundary as $\Sigma_0$ and satisfying $\xi\equiv -(c_o+t)$. Define 
$$\phi(t):=f(t)+th\:.$$

Note that 
$$0=\frac{d}{dt}\left[\Xi(t,f(t)\right]_{t=0}= P[h]+P[{\dot f}]+2=P[{\dot f}]\:,$$
where ${\dot f}=\partial_tf(t)_{t=0}$. Since ${\dot f}=0$ on $\partial \Sigma_0$, we have, by Lemma \ref{nas}, ${\dot f}\equiv 0$ and so $X_0+(th+f(t))\nu=X_0+th\nu+{\mathcal O}(t^2)$. {\bf q.e.d.}
\\

\begin{figure}[h!]
\makebox[\textwidth][c]{
\hspace{-1.1cm}
\begin{subfigure}[b]{0.45\linewidth}
\includegraphics[width=\linewidth]{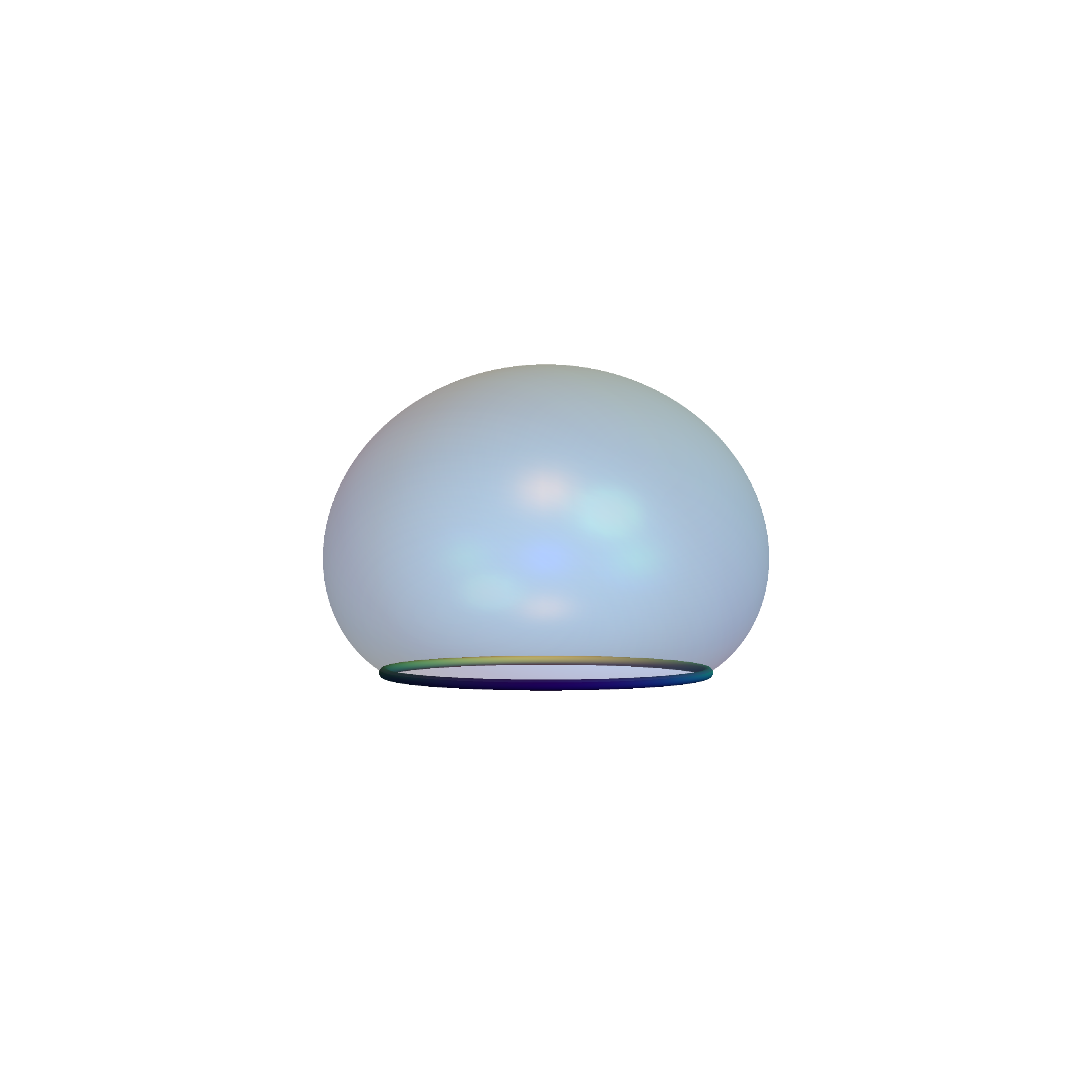}
\caption{$c_o=1.8$}
\end{subfigure}
\hspace{-1.6cm}
\begin{subfigure}[b]{0.45\linewidth}
\includegraphics[width=\linewidth]{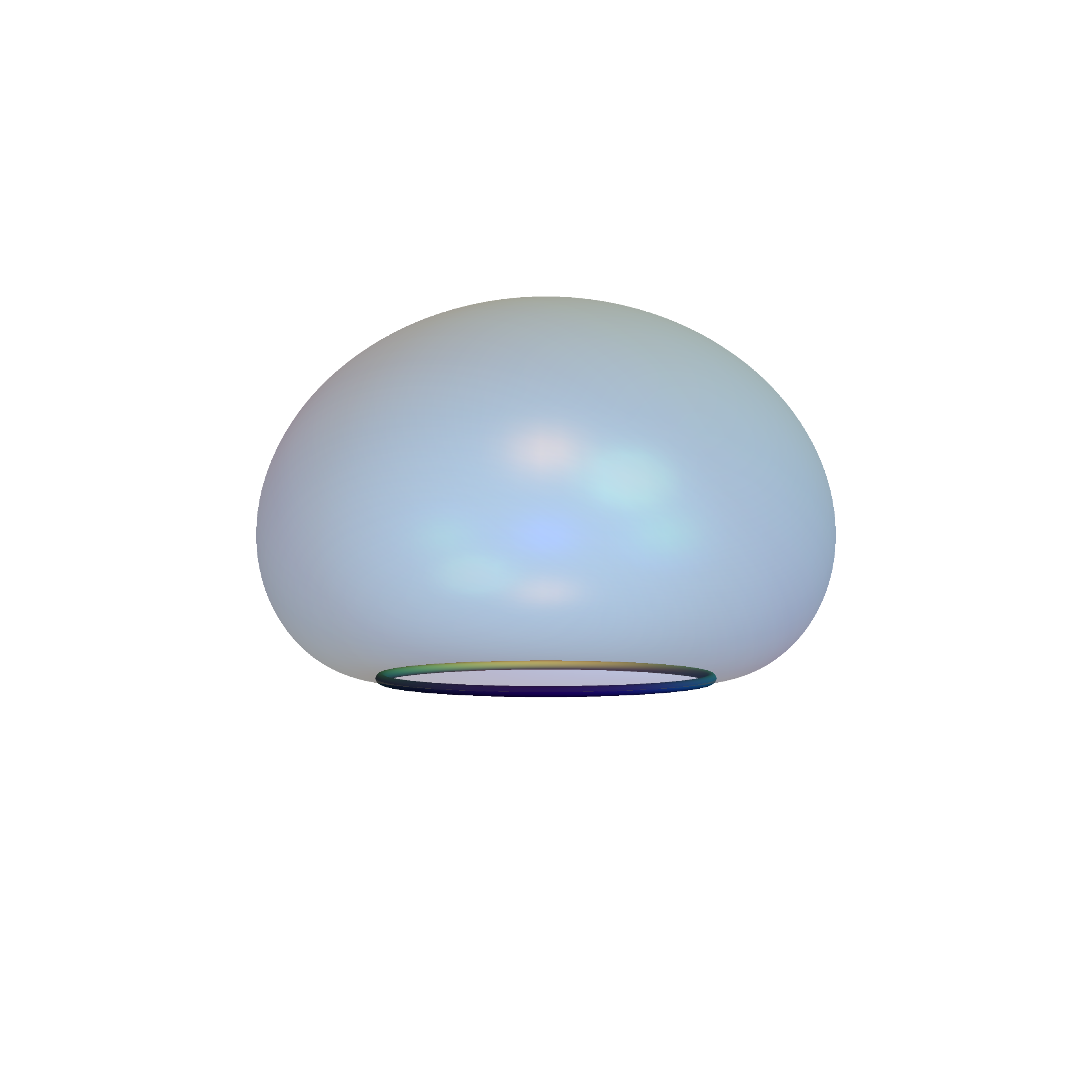}
\caption{$c_o=1.5$}
\end{subfigure}
\hspace{-0.9cm}
\begin{subfigure}[b]{0.45\linewidth}
\includegraphics[width=\linewidth]{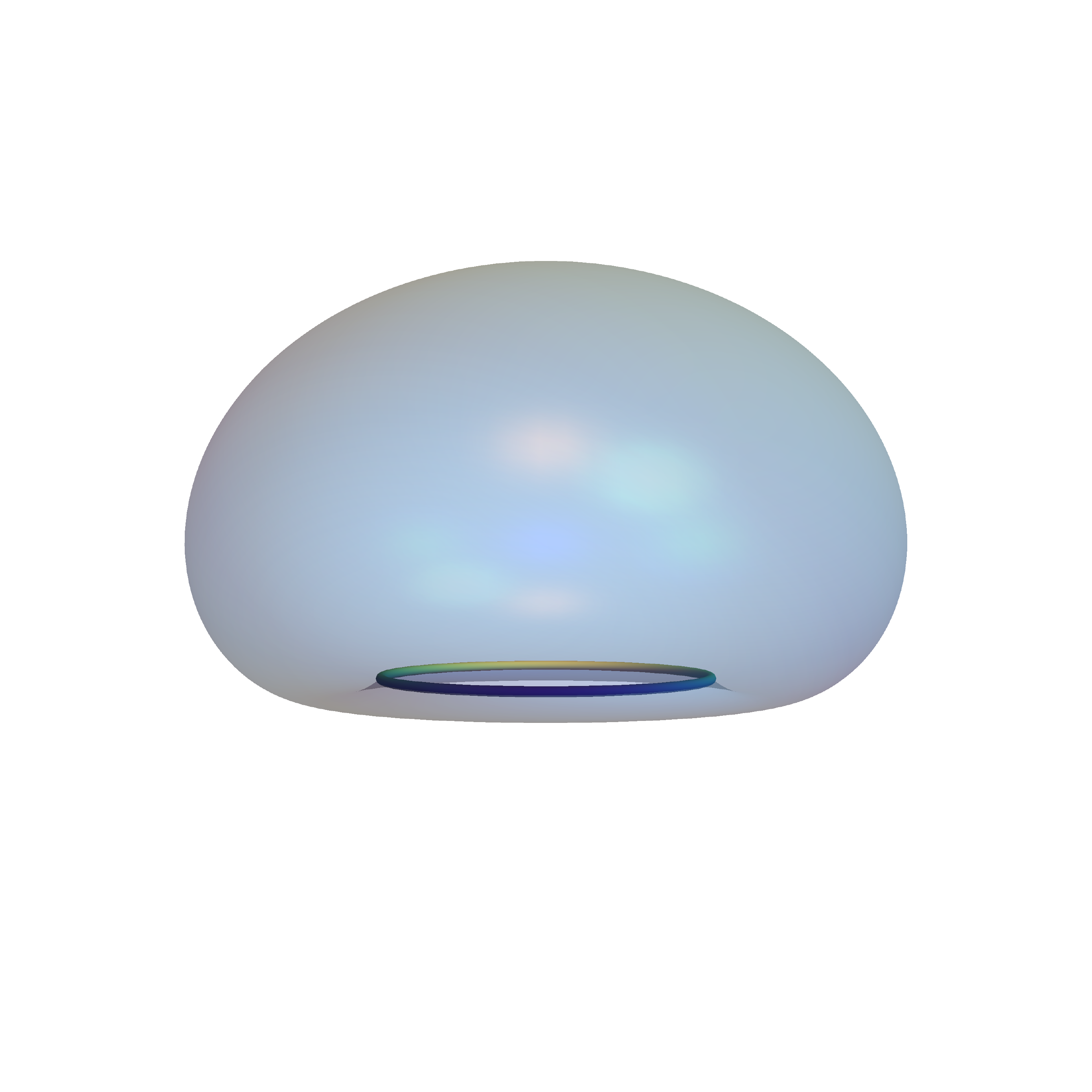}
\caption{$c_o=1.3$}
\end{subfigure}
\hspace{-0.4cm}
\begin{subfigure}[b]{0.45\linewidth}
\includegraphics[width=\linewidth]{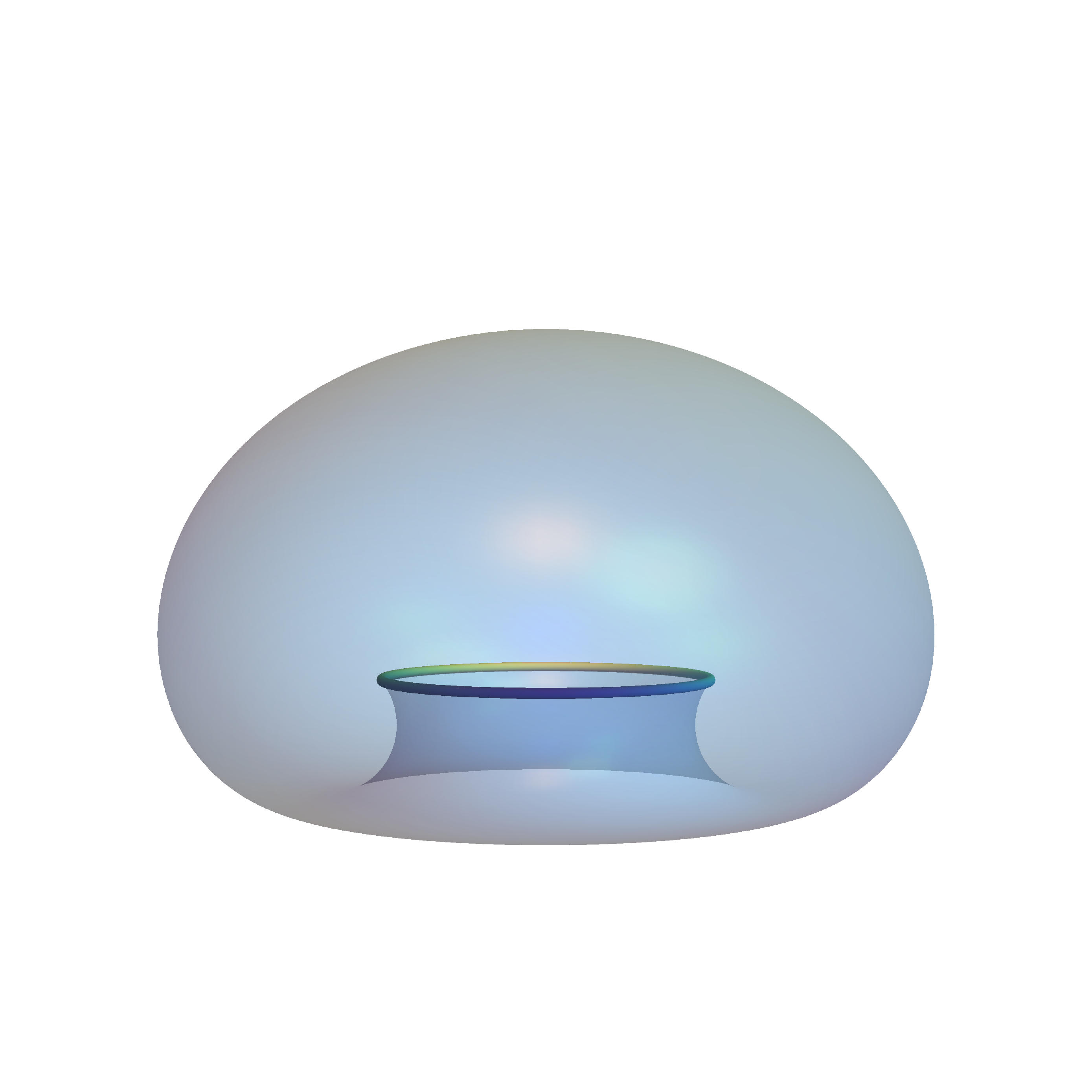}
\caption{$c_o=1.2$}
\end{subfigure}}
\caption{A family of axially symmetric surfaces $\{\Sigma_t\}$ whose boundary is the circle of radius $r=0.5$ located at $z=-3$. Here $\Sigma_0$ is the surface (B).}
\label{family}
\end{figure}

In Figure \ref{family} we show a family of axially symmetric surfaces $\Sigma_t$ satisfying \eqref{H+c} for different values of $c_o$ such that all of them share the same boundary circle. Observe that, as mentioned in the Introduction, the contact angle between the surface and the boundary circle varies.

\section{Producing the Bifurcation}

The basic abstract theory we will employ to produce a family of membranes bifurcating from the axially symmetric family $\{\Sigma_t\}$, $t\in(-\epsilon,\epsilon)$, is due to Crandall and Rabinowitz (\cite{CR}). We begin by introducing this result.

Let $Y$, $Z$ be real Banach spaces, $V$ an open neighborhood of $0$ in $Y$, $I$ a non-empty open interval, and $F:I\times V \to Z$ a twice continuously Fr\'echet differentiable mapping. For a linear mapping $T$, denote by $N[T]$ the kernel of $T$, and by $R[T]$ the image of $T$.

\begin{theorem}\label{CR}\mbox{\cite[Theorem~1.7]{CR}}\label{cran} Assume that $t_0 \in I$, $\eta\in V\subset Y$ and that the following statements hold:
\begin{enumerate}[(i)]
\item $F(t, 0)=0$ for all $t \in I$,
\item ${\rm dim} \left(N[D_\eta F(t_0, 0)]\right)={\rm codim} \left(R[D_\eta F(t_0, 0)]\right)=1$,
\item $D_t D_\eta F(t_0, 0)(y_0)\notin R[D_\eta F(t_0, 0)]$, where $y_0\in Y$ spans $N[D_\eta F(t_0, 0)]$.
\end{enumerate}
Let $W$ be any complement of ${\rm span}\,\{y_0\}$ in $Y$. Then there exists an open interval $\widehat{I}$ containing $0$ and continuously differentiable functions $t:\widehat{I}\to {\bf R}$ and $\mu:\widehat{I}\to W$ such that $t(0)=t_0$, $\mu(0)=0$, and if $y(s)=sy_0+s\mu(s)$, then $F(t(s), y(s))=0$. Moreover, $F^{-1}(\{0\})$ near $(t_0, 0)$ consists precisely of the curves $(t, 0)$, $t \in I$, and $(t(s), y(s))$, $s \in \widehat{I}$.
\end{theorem}

We now state and prove the main result of the paper, which shows the existence of the symmetry breaking bifurcation.
 
\begin{theorem}\label{bifurcation} Let $\Sigma_0$ be an axially symmetric convex embedded disc type surface satisfying
$$\xi=H+\frac{\nu_3}{z}=c_o\,,$$
and having a horizontal tangent plane along its boundary circle. Then there exists an open interval $\widehat{I}$ containing $0$ and a continuously differentiable function $t:\widehat{I}\to {\bf R}$ such that the one parameter family of non-axially symmetric surfaces $Y=Y(s):=X_0+\left(\phi(t(s))+s\nu_1+\mathcal{O}(s^2)\right)\nu$ satisfy $\xi\equiv\xi(t(s))$ and have the same boundary as $\Sigma_0$.

In some ${\mathcal C}^{2+\alpha}$ neighborhood of $\Sigma_0$, all surfaces with $\xi$ constant and having the same boundary as $\Sigma_0$ are equal, up to a rotation, to a surface in either the family $Y(s)$ or the family $X(t)$ of Theorem \ref{A}.
\end{theorem}

\begin{lemma}\label{hb} Let $\Sigma_0$ be as in Theorem \ref{bifurcation} and let the function $h(\varsigma)$, $\varsigma\in[0,\ell]$, be as in \eqref{h}. Then $h_\varsigma(0)\neq 0$.
\end{lemma}
{\it Proof.\:} We normalize $\psi$ by $\psi(0)=1$, then $h=c_o^{-1}\left(q(0)\psi-q\right)$. Both the inequalities $q(0)=X_0\cdot \nu(0)>0$ and $(X_0)_\varsigma\cdot X_0(0)>0$ are clear. Since $\kappa(0)<0$ holds (see Theorem \ref{shape}),
$$h_\varsigma(0)=c_o^{-1}\left(q(0)\psi_\varsigma(0)-q_\varsigma(0)\right)=c_o^{-1}\left(q(0)\psi_\varsigma(0)+\kappa(0)(X_0)_\varsigma\cdot X_0(0)\right)$$
will be negative if $\psi_\varsigma(0)\leq 0$ holds.

Note that since, $P[\nu_1]=0$,  $\nu_1\equiv 0$ on $\partial \Sigma_0$ and $\nu_1$ changes sign in $\Sigma_0$, the first Dirichlet eigenvalue of $P$ in $\Sigma_0$ is negative, so the solution $\psi$ of $P[\psi]=0$ must change sign in $\Sigma_0$. As shown in the proof of Lemma \ref{nas}, $\psi$ has no zeros above the circle in $\Sigma_0$ on which $\nu_3\equiv 0$ holds. Let $\widetilde{\varsigma}$ denote the smallest positive value of $\varsigma$ for which $\psi(\widetilde{\varsigma})=0$ holds so $\psi>0$ holds on $(0,\widetilde{\varsigma})$ and $\nu_3<0$ holds on $[0,\widetilde{\varsigma}]$. By \eqref{Pnu3} and $P[\psi]=0$, we then get
\begin{eqnarray*}
0&<&\int_0^{\widetilde{\varsigma}}\frac{-2\nu_3}{z^4}\:r\,d\varsigma=\int_0^{\widetilde{\varsigma}} \left(\psi \frac{P}{z^2}[\nu_3]-\nu_3  \frac{P}{z^2}[\psi]\right)r\,d\varsigma\\
&=& \int_0^{\widetilde{\varsigma}} \left(\psi\left[\frac{r(\nu_{3})_{\varsigma}}{z^2}\right]_\varsigma-\nu_3\left[\frac{r\psi_\varsigma}{z^2}\right]_\varsigma \right)r\,d\varsigma\\
&=& \int_0^{\widetilde{\varsigma}} \left(\frac{r}{z^2}\left[\psi (\nu_{3})_{\varsigma}-\nu_3\psi_\varsigma\right] \right)_\varsigma d\varsigma\\
&=& \left[\frac{r}{z^2}\left(\psi (\nu_{3})_{\varsigma}-\nu_3\psi_\varsigma\right)\right](\varsigma_1)-\left[\frac{r}{z^2}\left(\psi (\nu_{3})_{\varsigma}-\nu_3\psi_\varsigma\right)\right](0)\\
&=&\frac{r}{z^2}\left[-\nu_3\psi_\varsigma\right](\widetilde{\varsigma})-\left[\frac{r}{z^2}\left(\psi (\nu_{3})_{\varsigma}-\nu_3\psi_\varsigma\right)\right](0)\:.
\end{eqnarray*}
Note that $\psi_\varsigma(\widetilde{\varsigma})<0$ holds since $\psi$ changes from positive to negative at $\widetilde{\varsigma}$. Also $\left(\nu_{3}\right)_{\varsigma}(0)=-\kappa(0)z_\varsigma(0)=0$, so $\nu_3\psi_\varsigma(0)$ must be positive and hence $\psi_\varsigma(0)<0$ must hold. {\bf q.e.d.}
\\

{\it Proof of Theorem \ref{bifurcation}.\:} The procedure used to produce the bifurcation is  similar to that used in \cite{KPP} and \cite{P}. We take $V=\{\psi(r,\theta)\in\mathcal{C}_o^{2+\alpha}(\Sigma_0)\, \rvert\, \psi(r,-\theta)=\psi(r,\theta)\}$. For $Z$, we note that there is an embedding
$$\iota:\mathcal{C}_o^\alpha(\Sigma_0) \hookrightarrow (\mathcal{C}_o^\alpha(\Sigma_0))^\ast,$$
given by 
$$\langle \iota f, g\rangle:=\int_{\Sigma_0} f\,g\:z^{-2}\:d\Sigma\,,\quad\quad\quad f,g\in \mathcal{C}_o^\alpha (\Sigma_0)\:$$
and we let $Z=\iota(\mathcal{C}_o^\alpha(\Sigma_0))\subset (\mathcal{C}_o^\alpha(\Sigma_0))^\ast$. 

We recall from \cite{KPP} that for  functions $f$ in a sufficiently small neighborhood of $0\in \mathcal{C}^{2+\alpha}_o(\Sigma_0)$, every map $X_0+f\nu$ is an immersion. Conversely, if $\overline{X}$ is an immersion sufficiently close to $X_0$ in the $\mathcal{C}_o^{2+\alpha}$ topology and with $\overline{X}\equiv X_0$ on $\partial\Sigma_0$, then $\overline{X}=X_0+f\nu$ for some $f\in \mathcal{C}^{2+\alpha}_o(\Sigma_0)$.  
 
We define $F:I\times V \longrightarrow Z$ by
\begin{equation*}
F(t,\eta):=D_\eta\left({\widetilde{\mathcal{A}}}_{\phi(t)+\eta}-2c_o (t)\mathcal{U}_{\phi(t)+\eta}\right),
\end{equation*}
where $D_\eta$ denotes the derivative with respect to $\eta$ and the subscripts in the right hand side are used in place of the surface determined by the stated function. Here, the function $\phi(t)$ is the $\mathcal{C}^2$ curve defined in Theorem \ref{A}. From the first variation formula of \eqref{G}, computed in Theorem 4.2 of \cite{PP2}, we have for $f\in \mathcal{C}^\alpha_o(\Sigma_0)$
$$F(t,\eta)[f]=-2\int_{\Sigma_0} f\,\nu\cdot\nu_{\phi(t)+\eta}\left(\xi_{\phi(t)+\eta}+c_o(t)\right)z_{\phi(t)+\eta}^{-2}\:d\Sigma_{\phi(t)+\eta}\:,$$
so condition (i) of Theorem \ref{cran} is fulfilled since $\xi_{\phi(t)}\equiv -c_o(t)$.

For $\zeta\in V$, the derivative $D_\eta F(0, 0):V\longrightarrow Z$, is given as the action on $f\in \mathcal{C}^\alpha_o(\Sigma_0)$ by:
\begin{eqnarray*}
D_\eta F(0, 0)(\zeta)[f]&=&\partial_\rho \left[F(0,0+\rho \zeta)\right]_{\rho=0}[f]=-\int_{\Sigma_0} f\left(\delta_\zeta \xi\right)z^{-2}\,d\Sigma\nonumber\\
&=&-\int_{\Sigma_0} f\,P[\zeta]\:z^{-2}\,d\Sigma=-\int_{\Sigma_0} \zeta \,P[f]\:z^{-2}\,d\Sigma\:.
\end{eqnarray*}
Recall that the operator $z^{-2}P$ is self adjoint  and so  the last equality holds. By Lemmas \ref{nas} and \ref{m}, the kernel of $P$ is spanned by $\nu_1=z_\sigma \cos\theta$ and $\nu_2=z_\sigma \sin \theta$. However, the definition of $V$ restricts us to functions that are even in $\theta$ which makes the kernel one dimensional. Consequently, ${\rm dim}\left(N[D_\eta F(0, 0)]\right)=1$ since $\nu_1$ spans the null space. The statement about the codimension of the range follows from the Fredholm Alternative. This proves that condition (ii) of Theorem \ref{cran} is satisfied.

Next, we verify the condition (iii) of Theorem \ref{CR}. By the previous calculation
$$D_\eta F(t,0)(\zeta)[f]=-\int_{\Sigma_0}  f\,\nu_{\phi(t)}\cdot \nu\, P_{\phi(t)}\left[\zeta \nu_{\phi(t)}\cdot \nu\right]\:z_{\phi(t)}^{-2}\,d\Sigma_{\phi(t)}\,.$$
Therefore, using that $P_{\phi(0)}[\nu_1]=P[\nu_1]=0$, we get
\begin{equation}\label{Dth}
D_t D_\eta F(0,0)(\nu_1)[f]= -\int_{\Sigma_0}  f\,\partial_t\left(P_{\phi(t)}\right)_{t=0}[ \nu_1]\:z^{-2}\,d\Sigma\:.
\end{equation}

Note that $P_{\phi(t)}[E_1\cdot \nu_{\phi(t)}] \equiv 0$, so
\begin{equation}\label{P'}
\partial_t\left(P_{\phi(t)}\right)_{t=0}[\nu_1]+P[ \nu_1']=0\:,
\end{equation}
where $\nu_1'=\partial_t\left(E_1\cdot \nu_{\phi(t)}\right)_{t=0}$.

Suppose to the contrary that  $D_tD_\eta F(0,0)(\nu_1)$ is in the range of $D_\eta F(0,0)$. Then, by \eqref{Dth}, there exists $\zeta\in V$ such that for all $\psi\in V$
$$\int_{\Sigma_0}  \psi\,\partial_t\left(P_{\phi(t)}\right)_{t=0}[\nu_1]\:z^{-2}\,d\Sigma=\int_{\Sigma_0} \psi\,P[\zeta]\:z^{-2}\,d\Sigma\:,$$
holds. From this we get
$$D_tD_\eta F(0,0)(\nu_1)[\nu_1]=\int_{\Sigma_0} \nu_1\,P[\zeta]\:z^{-2}\,d\Sigma=\int_{\Sigma_0}\zeta\,P[\nu_1]\:z^{-2}\,d\Sigma=0\,.$$
On the other hand, by \eqref{P'}, we get
\begin{eqnarray*}
D_tD_\eta F(0,0)(\nu_1)[\nu_1]&=&\int_{\Sigma_0} \nu_1\partial_t\left(P_{\phi(t)}\right)_{t=0}[\nu_1]\:z^{-2}\,d\Sigma \\
&=&-\int_{\Sigma_0} \nu_1\,P[\nu_1']\:z^{-2}\,d\Sigma\\
&=&\oint_{\partial \Sigma_0} \nu_1'\,\partial_n \nu_1\:z^{-2}\,ds\:.
\end{eqnarray*}

The proof of Theorem \ref{A} implies that $\nu_1'=-\nabla h\cdot  \nabla x=- h_\varsigma r_\varsigma \cos \theta $, while $\partial_n \nu_1=-z_{\varsigma \varsigma} \cos \theta=\kappa r_\varsigma \cos \theta$. The integrand in the last integral above is therefore a non zero multiple of $h_\varsigma(0) \cos^2\theta$ so a contradiction is reached due to Lemma \ref{hb}. This verifies condition (iii) of Theorem \ref{CR} and finishes the proof of the statement.
{\bf q.e.d.}
\\

\begin{figure}[h!]
\makebox[\textwidth][c]{
\begin{subfigure}[b]{0.37\linewidth}
\includegraphics[width=\linewidth]{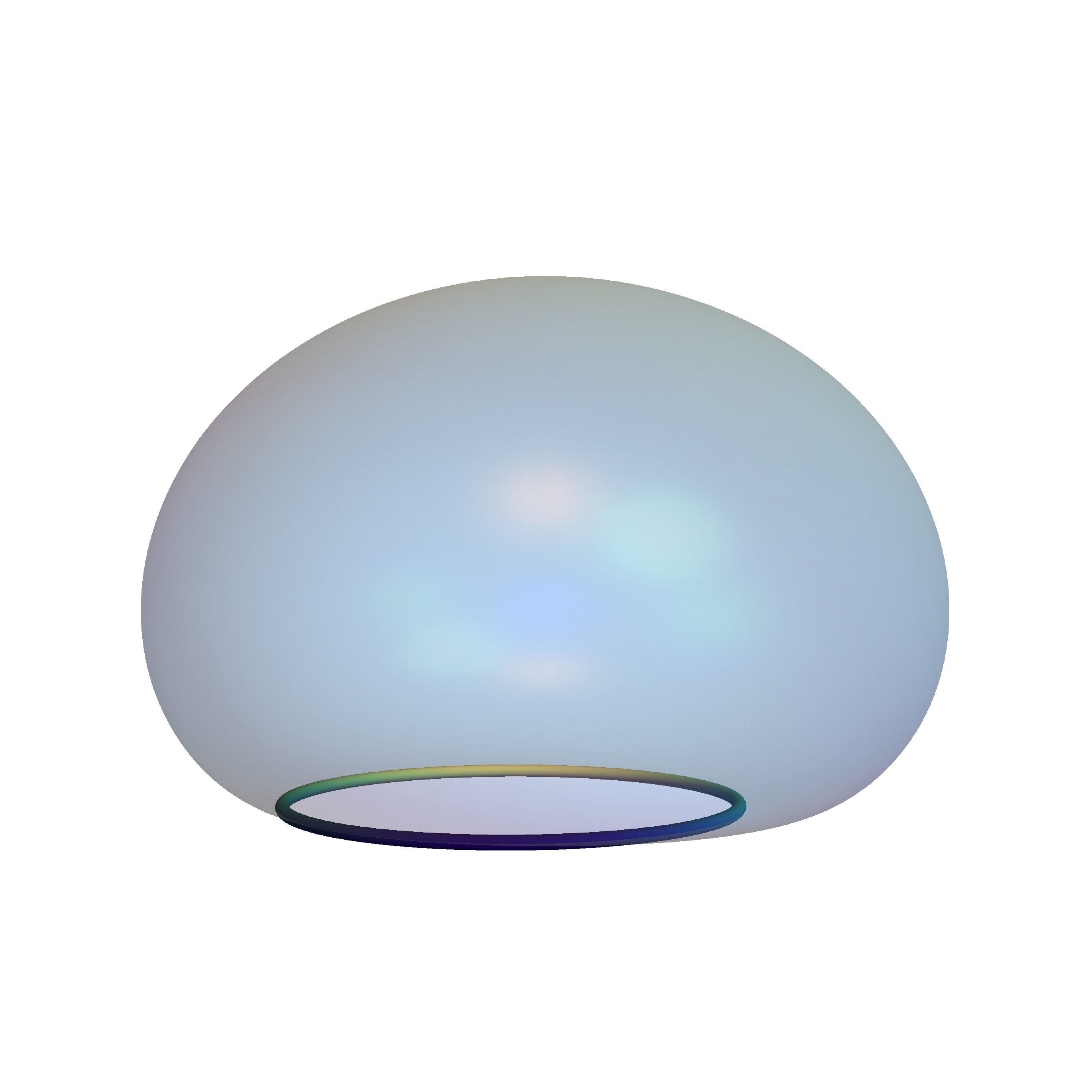}
\end{subfigure}
\begin{subfigure}[b]{0.37\linewidth}
\includegraphics[width=\linewidth]{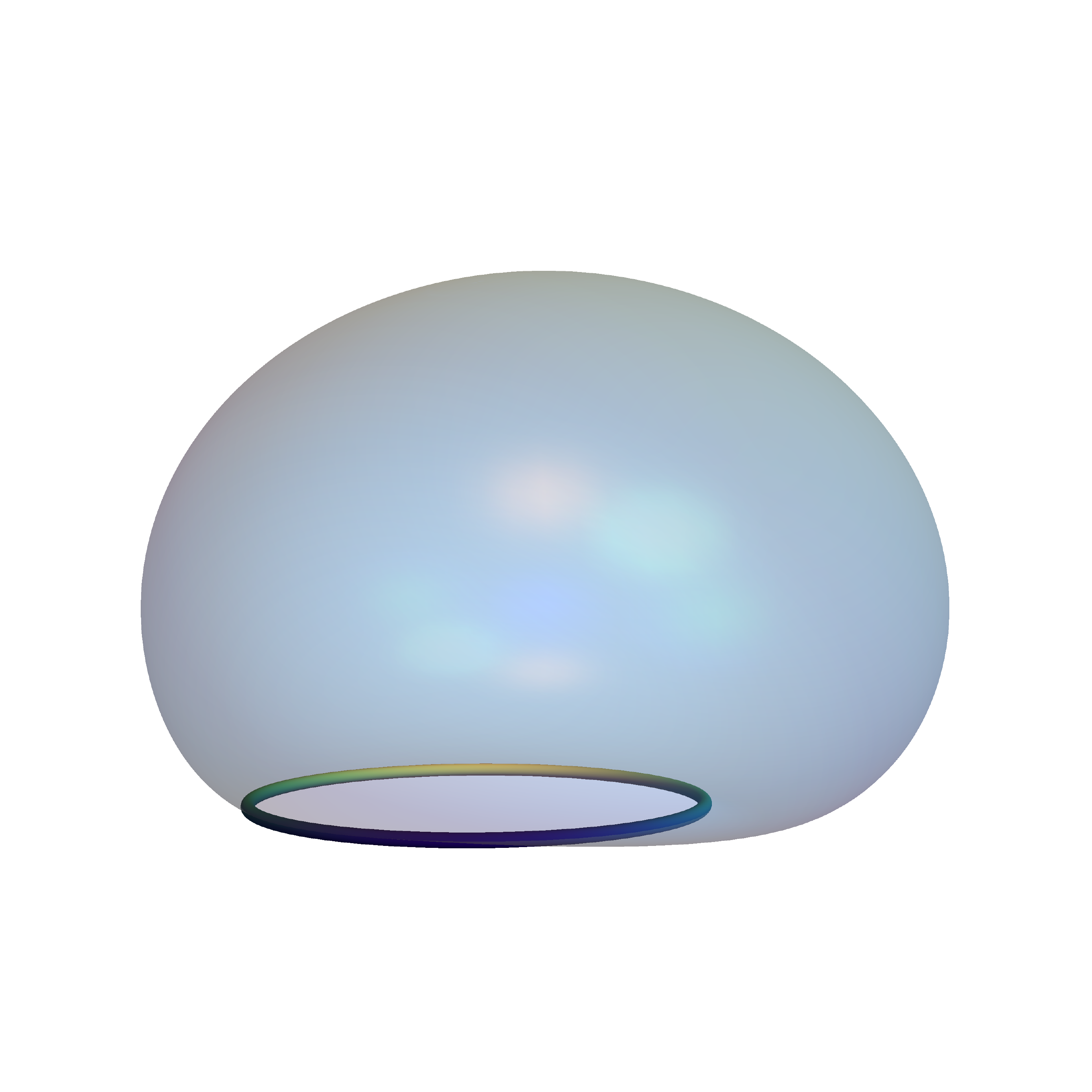}
\end{subfigure}
\begin{subfigure}[b]{0.37\linewidth}
\includegraphics[width=\linewidth]{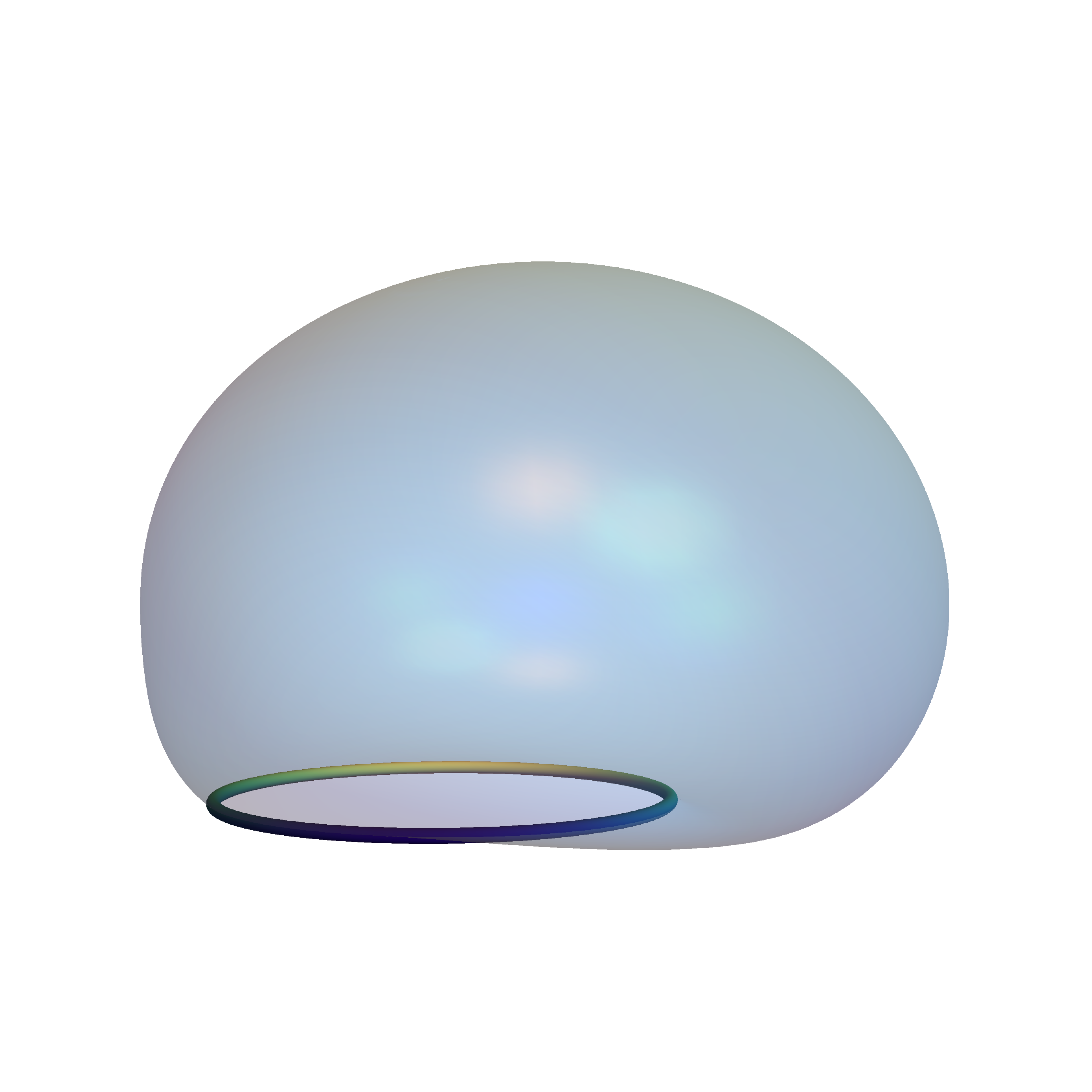}
\end{subfigure}
\begin{subfigure}[b]{0.37\linewidth}
\includegraphics[width=\linewidth]{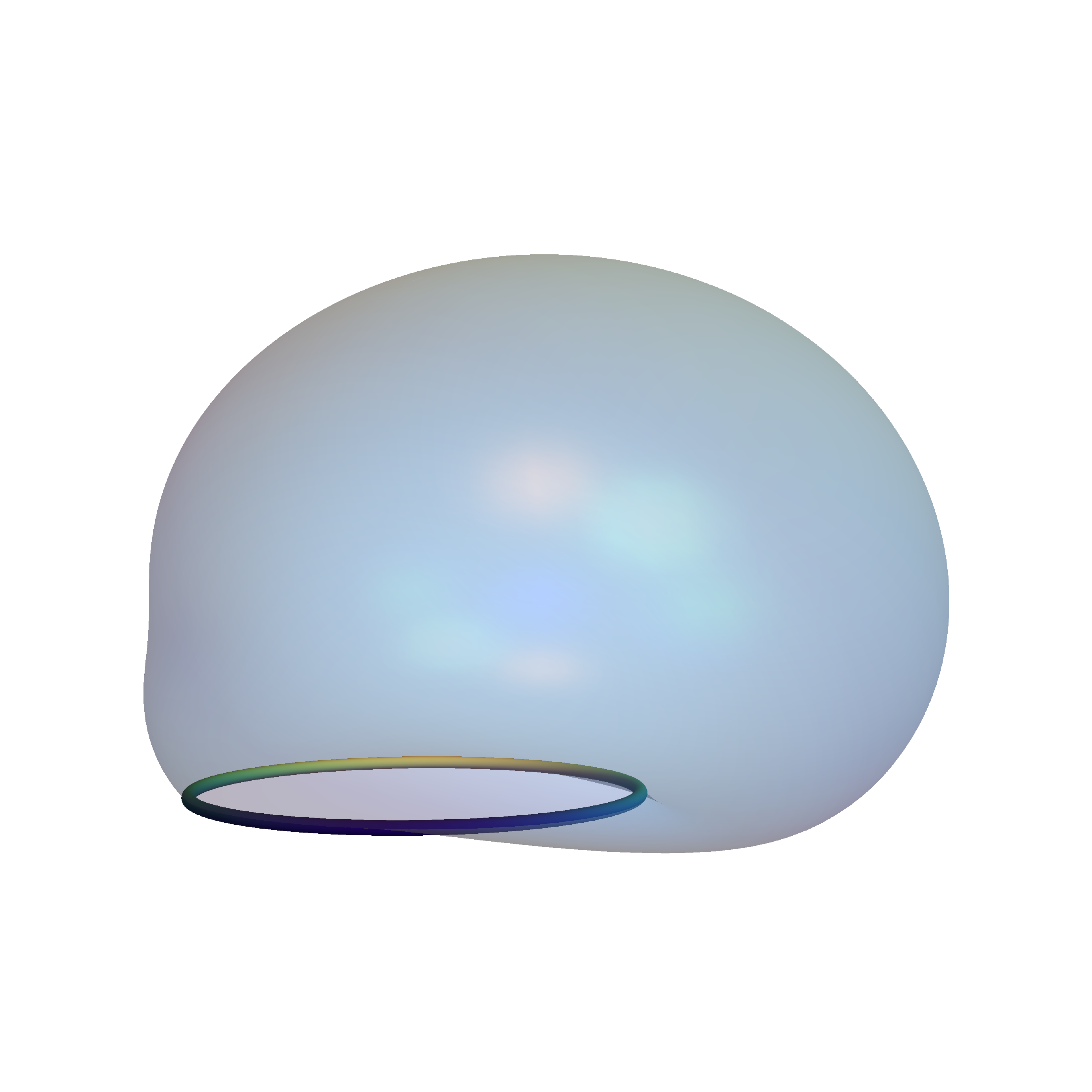}
\end{subfigure}}
\caption{Linear approximations of surfaces in the bifurcating branch. The corresponding bifurcation surface $\Sigma_0$ is shown in Figure \ref{family}, (B).}
\label{NewFamily}
\end{figure}

In Figure \ref{NewFamily} we illustrate the linear approximation of the non-axially symmetric surfaces satisfying \eqref{H+c} which bifurcate from the family of Figure \ref{family} at $\Sigma_0$ (Figure \ref{family}, (B)).

We will next numerically examine an interesting behavior of $h_\varsigma(0)\neq 0$ for various examples. To do this we modify the system \eqref{odesys} to contain the function $h$ and its first derivative
\begin{eqnarray*}\label{odesys2}
r_\varsigma&=&\cos\varphi\,,\nonumber\\
z_\varsigma&=&\sin\varphi\,,\\
\varphi_\varsigma&=&-2\frac{\cos\varphi}{z}-\frac{\sin\varphi}{r}+2c_o\:,\\
h_\varsigma&=&w\\
w_\varsigma&=&-\left(\lVert d\nu\rVert^2-2\left[\frac{\cos\varphi}{z}\right]^2\right)h-\left(\frac{\cos\varphi}{z}-2\frac{\sin \varphi}{r}\right)w-2\:,
\end{eqnarray*}
where
$$\lVert d\nu\rVert^2=4H^2-2K=\frac{\sin^2 \varphi}{r^2}+\left(-2\frac{\cos\varphi}{z}-\frac{\sin\varphi}{r}+2c_o\right)^2.$$
The last two equations in the system above represent the second order ordinary differential equation $P[h]=-2$.

The values of $h_\varsigma(0)$ are given in Table \ref{tab2} while the corresponding profile curves are shown in Figure \ref{Profiles}. The data indicates that as the surface becomes more spherical, the values of $h_\varsigma(0)$ tend to zero through negative values.

\begin{center}
\begin{table}[H]
\caption{Numerical values of $h_\varsigma(0)$ for the generating curves illustrated in Figure \ref{Profiles}.} \label{tab2}
\makebox[\textwidth][c]{
\def\arraystretch{2}
\begin{tabular}{| c | c | c | c | c | c |} 
\hline
$z_o=z(\ell)$ & $-0.55$ (Blue) & $-0.6$ (Green) & $-0.7$ (Yellow) & $-0.9$ (Red) & $-1.2$ (Purple) \\
\hline
$h_\varsigma(0)$ & $-23.1896$ & $-13.577$ & $-7.3487$ & $-3.8685$ & $-2.3639$ \\
\hline
\end{tabular}}
\end{table}
\end{center}

\section*{Acknowledgments}

The second author has been partially supported by the AMS-Simons Travel Grants Program 2021-2022. He would also like to thank the Department of Mathematics and Statistics of Idaho State University for its warm hospitality.

\bigskip

\begin{flushleft}
Bennett P{\footnotesize ALMER}\\
Department of Mathematics,
Idaho State University,
Pocatello, ID 83209,
U.S.A.\\
E-mail: palmbenn@isu.edu
\end{flushleft}

\bigskip

\begin{flushleft}
\'Alvaro P{\footnotesize \'AMPANO}\\
Department of Mathematics and Statistics, Texas Tech University, Lubbock, TX
79409, U.S.A.\\
E-mail: alvaro.pampano@ttu.edu
\end{flushleft} 

\end{document}